\documentclass[12pt,reqno]{article}

\usepackage[usenames]{color}
\usepackage[dvipsnames]{xcolor}
\usepackage{amssymb}
\usepackage{amsmath}
\usepackage{amsthm}
\usepackage{amsfonts}
\usepackage{amscd}
\usepackage{graphicx}
\usepackage{algorithm}
\usepackage{algorithmic}

\usepackage[colorlinks=true,
linkcolor=webgreen,
filecolor=webbrown,
citecolor=webgreen]{hyperref}

\definecolor{webgreen}{rgb}{0,.5,0}
\definecolor{webbrown}{rgb}{.6,0,0}

\usepackage{color}
\usepackage{fullpage}
\usepackage{float}

\usepackage{graphics}
\usepackage{latexsym}
\usepackage{epsf}

\newcommand{\seqnum}[1]{\href{https://oeis.org/#1}{\underline{#1}}}

\usepackage{enumerate,tikz,listings,array,longtable}
\usetikzlibrary{shapes,arrows,chains}
\usepackage{pgf}

\DeclareMathOperator{\IMA}{im}
\DeclareMathOperator{\Image}{Im}
\DeclareMathOperator{\FX}{Fx}
\DeclareMathOperator{\fx}{fx}
\DeclareMathOperator{\RmL}{RmL}
\DeclareMathOperator{\mb}{mb}
\DeclareMathOperator{\NMB}{Nmb}
\DeclareMathOperator{\nmb}{nmb}
\DeclareMathOperator{\MB}{Mb}
\DeclareMathOperator{\sueq}{suc}
\DeclareMathOperator{\Succeq}{Suc}
\DeclareMathOperator{\SF}{SF}
\DeclareMathOperator{\SP}{SP}
\DeclareMathOperator{\BP}{BP}
\DeclareMathOperator{\blocks}{bl}
\DeclareMathOperator{\AX}{Ax}
\DeclareMathOperator{\ax}{ax}
\DeclareMathOperator{\ExcL}{ExcL}
\DeclareMathOperator{\AXL}{AxL}
\DeclareMathOperator{\Exc}{Exc}
\DeclareMathOperator{\exc}{exc}
\DeclareMathOperator{\MX}{Mx}
\DeclareMathOperator{\wExc}{Wex}
\DeclareMathOperator{\wExcL}{WexL}
\DeclareMathOperator{\wexc}{wex}
\DeclareMathOperator{\RGF}{RGF}
\DeclareMathOperator{\cyc}{cyc}
\DeclareMathOperator{\Swap}{Swap}

\usepackage{datetime}
\newdate{date}{15}{03}{2022}
\date{\displaydate{date}}
\begin{document}
	
	\theoremstyle{plain}
	\newtheorem{theorem}{Theorem}
	\newtheorem{corollary}[theorem]{Corollary}
	\newtheorem{lemma}[theorem]{Lemma}
	\newtheorem{proposition}[theorem]{Proposition}
	
	\theoremstyle{definition}
	\newtheorem{definition}[theorem]{Definition}
	\newtheorem{example}[theorem]{Example}
	\newtheorem{conjecture}[theorem]{Conjecture}
	\newtheorem{problem}[theorem]{Problem}
	\newtheorem{remark}[theorem]{Remark}
	\newtheorem{question}[theorem]{Question}
	
	\begin{center}
		{\large \bf Set Partitions and Other Bell Number Enumerated Objects}\\
		
		\bigskip
		Fufa Beyene\\
		Department of Mathematics, Addis Ababa University\\
		1176 Addis Ababa, Ethiopia\\
		\href{mailto: fufa.beyene@aau.edu.et}{\tt fufa.beyene@aau.edu.et}\\
		\ \\
		J\"orgen Backelin\\
		Department of Mathematics, Stockholm University\\
		SE-106 91 Stockholm, Sweden\\
		\href{mailto: joeb@math.su.se}{\tt joeb@math.su.se}\\
		\ \\
		Roberto Mantaci	\\
		IRIF, Universit\'e de Paris\\8 Place Aur\'elie Nemours\\ F-75013 Paris, France\\
		\href{mailto: mantaci@irif.fr}{\tt mantaci@irif.fr}\\
		\ \\
		Samuel A. Fufa\\
		Department of Mathematics, Addis Ababa University\\
		1176 Addis Ababa, Ethiopia\\
		\href{mailto: samuel.asefa@aau.edu.et}{\tt samuel.asefa@aau.edu.et}\\
		%\today
	\end{center}
	\begin{abstract}
		In this paper, we study classes of subexcedant functions enumerated by the Bell numbers and present bijections on set partitions. We present a set of permutations whose transposition arrays are the restricted growth functions, thus defining Bell permutations of the second kind. We describe a bijection between Bell permutations of the first kind (introduced by Ponti and Vajnovzski) and the second kind. We present two other Bell number enumerated classes of subexcedant functions. Further, we present bijections on set partitions, in particular, an involution that interchanges the set of merging blocks and the set of successions. We use the bijections to enumerate the distribution of these statistics over the set of set partitions, and also give some enumeration results.
	\end{abstract}
	
	\section{Introduction}
	Let $n$ be a fixed positive integer and let $[n]:=\{1, 2, \ldots, n\}$.  A set partition of $[n]$ is a collection of pairwise disjoint non-empty subsets of $[n]$ such that their union forms the whole set $[n]$.
	
	For any set $S$ the function $\sigma : [n]\longrightarrow S$ corresponds to the word $\sigma(1)\sigma(2)\cdots\sigma(n)$. In particular, a permutation is a word with distinct symbols. 
	
	Permutations and set partitions are among the richest objects in enumerative combinatorics. A basic reason for this fact is the wide variety of ways to represent a permutation and a set partition combinatorially.
	A second reason for their richness is the wide variety of interesting \emph{statistics}. Among many statistics on the set of permutations, the most classical ones are Eulerian and Mahonian statistics. Some of these are descents, weak excedances, anti-excedances, right-to-left minima or maxima, inversions, etc. On the other hand, we recall the two most basic enumerations for set partitions: the total number of set partitions over $[n]$ and the number of set partitions over $[n]$ having $k$ blocks are the Bell number, $B(n)$ and the Stirling number of the second kind, $S(n,k)$ respectively (see \cite{Bo,Ma,Ro,St1}).
	
	Both permutations and set partitions can be coded by subexcedant functions, i.\ e., functions $f: [n]\mapsto[n]$ such that  $1\leq f(i)\leq i$ for all $i \in [n]$ (in some contexts it is rather required  that $0\leq f(i)\leq i-1$).
	
	Some permutation codes with subexcedant functions are very well known (Lehmer code or inversion table, Denert code, and so on \cite{Le,Du-Vi,Fo-Ze,Ma-Ra}). On the other hand, a way to code set partitions with subexcedant functions is given by Mansour's definition of canonical form for a set partition $P$ in the standard form, the elements in each block are arranged increasingly, and the blocks are arranged in increasing order of their minima (see \cite{Ma}). In the canonical form, any integer $i\in [n]$ is coded with the index of the block of $P$ where it belongs, where $P$ is in its standard form. In fact, canonical forms of set partitions are restricted growth functions ($\RGF$).
	
	Several  properties of set partitions or permutations can be read easily from their corresponding codes, which allows one to prove some results elegantly by reasoning on the codes rather than the coded objects. See for instance, the article of Baril and Vajnovszki \cite{Ba-Va}, and also the article of Foata and Zeilberger \cite{Fo-Ze}.
	
	Mantaci and Rakotondrajao (see \cite{Ma-Rak})  studied the bijection $\phi$ associating a subexcedant function $f$ with the permutation $\sigma =\phi(f)= (n,f_n)(n{-}1,f_{n-1})\cdots(1,f_1)$, where $f_i=f(i),~ \forall i\in[n]$ and related the image values of $f$ to the anti-excedances of $\sigma$. Later, Baril \cite{Ba1} independently studied a variation of the bijection $\phi$, here denoted by $\chi$,	given by simply inverting the order of the product of transpositions in the definition of $\phi$, and he called the subexcedant function associated with a permutation via this bijection the \emph{transposition array}. Baril \cite{Ba} also studied, in particular, the positions of weak excedances in a permutation using the corresponding subexcedant function.
	
	Mansour and Munagi \cite{Ma-Mu} studied set partitions according to the number of circular successions, i.\ e., the number of consecutive element pairs inside a block assuming that the elements are arranged around a circle. Callan \cite{Ca} has proved that the statistics for the number of singletons in all set partitions is equal to the one for the number of circular successions, by giving a bijection in terms of an algorithm that interchanges singletons and circular successions. Callan also proved that his bijection is an involution on set partitions and that it preserves the non-crossing partitions.
	
	In this paper, we study families of subexcedant functions that are enumerated by the Bell numbers. We present bijections between these classes and set partitions. We enumerate these classes based on certain statistics. We also present an involution on set partitions and use it to give some enumeration results. 
	
	In Section \ref{bellperm}, we study a class of permutations whose transposition array is the restricted growth functions, we will call this class Bell permutations of the second kind. We prove that the statistic of the number of weak excedances is the Stirling number of the second kind, enumerate the statistic of the number of cycles, and their joint distribution. In Section \ref{otherclassSef}, we present two more families of subxcedant functions enumerated by the Bell numbers and bijections between these classes and set partitions.
	
	In Section \ref{bijBellperm1and2}, we provide a bijection between permutations of the second kind and another Bell-counted class of permutations introduced by Poneti and Vajnovszki \cite{Po-Va}. 
	
	Finally, in Section \ref{involution}, we present bijections on set partitions, in particular we present an involution that interchanges the number of merging blocks and the number of successions. We use the bijections to give some enumeration results. And also the generating function for the joint distribution of these statistics.
	
	\section{Notation and Preliminaries}
	\subsection*{Permutations}
	Recall that a permutation over $[n]$ is a bijection $\sigma:[n]\mapsto[n]$. Let $\mathfrak{S}_n$ denote the set of all permutations over $[n]$. A permutation $\sigma\in\mathfrak{S}_n$ can be written as a word $\sigma=\sigma(1)\sigma(2)\ldots\sigma(n)$ (whence the $\sigma(i)$ also are called \emph{letters}), or in cycle notation as a product of disjoint cycles,  where as usual a \emph{cycle} in $\sigma$ can be written as $(j,\sigma(j),\sigma^2(j), \ldots, \sigma^{t-1}(j))$, where $t$, the \emph{length} of the cycle, is the smallest positive integer such that $\sigma^t(j)=j$. Cycles of length one are \emph{fixed points}. The cycle notation is noted $\sigma=C_1C_2\cdots C_k$, where the $C_i$'s are disjoint cycles and the minima of the cycles form an increasing sequence. We let $\cyc(\sigma)$ denote the number of cycles of $\sigma$. A \emph{transposition} is a permutation that swaps two integers and fixes all the others.
	
	We say that a permutation $\sigma$ over $[n]$ has an \textit{excedance (weak excedance, anti-excedance)} in a position $i$ if $\sigma(i)>i$ ($\sigma(i)\ge i, \sigma(i)\le i$, respectively), where $i\in[n]$. We use the notation 
	\begin{align*} &\Exc(\sigma):=\{i : 1\leq i\leq n, \sigma(i)>i\},\\
		&\wExc(\sigma):=\{i : 1\leq i\leq n, \sigma(i)\geq i\}, \text{ and }\\
		&\AX(\sigma):=\{i : 1\leq i\leq n, \sigma(i)\leq i\}.
	\end{align*}
	We also use the notation 
	$\exc(\sigma):=|\Exc(\sigma)|$, $\wexc(\sigma):=|\wExc(\sigma)|,$ and 
	$\ax(\sigma):=|\AX(\sigma)|.$
	
	The set of \textit{excedance letters} (\textit{weak excedance letters, anti-excedance letters}) of $\sigma$ is defined as $\ExcL(\sigma):=\{\sigma(i): i\in\Exc(\sigma)\}$ ($\wExcL(\sigma):=\{\sigma(i): i\in\wExc(\sigma)\}$, $\AXL(\sigma):=\{\sigma(i): i\in\AX(\sigma)\}$, respectively).
	
	\subsection*{Subexcedant functions}
	We let $\SF(n)$ denote the set of all subexcedant functions over $[n]$. For $f=f_1f_2\cdots f_n\in \SF(n)$, we use the notation: $\Image(f):=f([n])$, the image set of $f$, and $\IMA(f):=|\Image(f)|$. We say that $f$ has a \emph{leftmost (rightmost) occurrence} in a position $i$ if $f_i\notin\{f_1, \ldots, f_{i-1}\}$, i.\ e., $i=\min(f^{-1}(f_i))$ (or $f_i\notin\{f_{i+1}, \ldots, f_n\}$, i.\ e., $i=\max(f^{-1}(f_i))$, respectively), where $i\in[n]$. If $i$ is a leftmost (rightmost) occurrence in $f$, then we say that $f_i$ is a leftmost (or rightmost) letter. The set of fixed points of $f$ is given by
	$$\FX(f):=\{i : 1\leq i\leq n, f_i=i\}.$$ We let $\fx(f):=|\FX(f)|$.
	\subsection*{Set Partitions}
	A set partition $P$ of $[n]$ is defined as a collection $B_1, \ldots, B_k$ of nonempty disjoint subsets of $[n]$ such that $\bigcup_{i=1}^kB_i=[n]$. The subsets $B_i$ will be referred to as \emph{blocks}. The \emph{block representation} $P=B_1|B_2|\cdots|B_k$ of a set partition $P$ is said to be \emph{standard} if the blocks $B_1, \ldots, B_k$ are sorted in such a way that $\min(B_1) < \min(B_2) <\cdots < \min(B_k)$ and if the elements of each block are arranged in increasing order. 
	
	We consider set partitions only in their standard representation. 
	
	We let $\SP(n)$ denote the set of all set partitions over $[n]$. We also let $\blocks(P)$ denote the number of blocks of a set partition $P$ and $\SP(n, k):=\{P\in\SP(n): \blocks(P)=k\}$. 
	
	Recall that $|\SP(n, k)|=S(n,k)$, where $S(n, k)$ is the Stirling number of the second kind.
	
	For $2\leq i\leq k$, we say that the block $B_i$ is \emph{merging} if $\max(B_{i-1})<\min(B_i)$. A set partition without merging blocks is called \emph{merging-free}. 
	
	If the integers of the pair $({a-1},a)$, where $a>1$, are in the same block of $P$, then $a$ is said to be a \emph{succession} of $P$. In literature, the name  ``succession" is used for the first element of the pair, but for our purposes, we prefer to use it for the second element. 
	
	We let $\MB(P)$, $\Succeq(P)$, and $\NMB(P)$ denote the set of the minimum elements of merging blocks, the set of successions of $P$, and the set of the minimum elements of non-merging blocks, respectively. We use the notation $\mb(P):=|\MB(P)|,~ \sueq(P):=|\Succeq(P)|$, and $\nmb(P):=|\NMB(P)|$.
	\begin{remark}	\label{eq0}
		For any $P\in\SP(n)$, every element $i$ of $[n]$ is necessarily in one of the first $i$ blocks of $P$. 
	\end{remark}
	The \emph{canonical form} of a set partition $P=B_1|B_2|\cdots|B_k$ is an $n$-tuple $f=f_1f_2 \cdots f_n$ indicating for each integer $j$ the index of the block in which it occurs, i.\ e., $B_j= f^{-1}(j)$ for all $j\in[k]$. For instance, the canonical form of $P=1~5~7|2~4|3~8|6\in\SP(8)$ is $f=12321413$.
	
	\begin{remark}
		The block $B_i$ contains its own index $i$ if and only if $i\in\FX(f)$. 
	\end{remark}
	Note that the canonical form of a set partition is a subexcedant function, but not all subexcedant functions are canonical forms of set partitions.
	
	A \emph{restricted growth function} $(\RGF)$ over $[n]$ is a function $f :[n]\mapsto[n]$, where $f=f_1\cdots f_n$ such that $f_1=1$ and $f_i \leq 1+\max\{f_1,\ldots, f_{i-1}\}$ for $2\leq i\leq n$, or equivalently, such that the set $\{f_1, f_2, \ldots, f_i\}$ is an integer interval for all $i\in[n]$. The canonical forms of set partitions are exactly the restricted growth functions ($\RGF$) (see \cite[p. 2]{Ma-Ra}). We let $\RGF(n)$ denote the set of all restricted growth functions over $[n]$.
	\section{Bell Permutations of the second kind}
	\label{bellperm}
	In this section, we study the class of permutations associated to $\RGF$s under $\chi$, the bijection given by Baril \cite{Ba1}. This set of permutations is counted by the Bell numbers, therefore we will call these objects ``Bell permutations of the second kind'' (Poneti and Vajnovszki in \cite{Po-Va} already introduced  another family of permutations counted by the Bell numbers that they called ``Bell permutations''). 
	
	The bijection $\chi$ is given by $\chi: \SF(n)\mapsto \mathfrak{S}_n$, where the permutation $\sigma=\chi(f)$ is defined by the product of transpositions:
	$$\sigma=(1,f_1)(2,f_2)\cdots(n,f_n),$$
	where the product is taken from right-to-left.
	The subexcedant function $f=\chi^{-1}(\sigma)$ is called the transposition array of $\sigma$. It is shown in \cite{Ba} that $\Image(f)=\wExc(\sigma)$. For instance, take $f=121132342\in \SF(9)$. Then \begin{align*}\sigma=\chi(f)&=(1,1)(2,2)(3,1)(4,1)(5,3)(6,2)(7,3)(8,4)(9,2)\\&=497812536,
	\end{align*}
	and $\Image(f)=\{1,2,3,4\}=\wExc(\sigma)$.
	\begin{remark}(\cite{Ba})\label{propertyofphitilde1}  
		The rightmost occurrences of $f$ are the weak excedance letters of $\chi(f)$.
	\end{remark}
	\begin{remark}\label{propertyofphitilde2}
		Let $f=f_1f_2\cdots f_n$ and $\sigma=\chi(f)$. We have $i\in\FX(f)$ if and only if $i$ is the minimum element of some cycle of $\sigma$.
	\end{remark}
	In \cite{Be-Ma2}, the following was essentially proved:
	\begin{lemma}\label{lemmaInom}
		Let $\sigma=\sigma(1)\sigma(2)\cdots\sigma(n)\in\mathfrak{S}_n$. If $f=\chi^{-1}(\sigma)=f_1f_2...f_n$, then $f(i)~=~\sigma^{-t}(i)\leq~i$, where $t\geq1$ is chosen as small as possible.\qed
	\end{lemma}
	The following proposition presents an alternative algorithm to compute $\sigma$ from $f$ as a product of disjoint cycles.
	\begin{proposition}\label{propcycsigma}
		If $f\in \SF(n)$,  then $\sigma =\chi(f)$ can be constructed as follows. For $i~=~1,2,\ldots,n$:
		\begin{itemize}
			\item if $f_i = i$, then add a new singleton cycle: $(i)$,
			\item if $f_i < i$, then insert $i$ after $f_i$ in its cycle.
		\end{itemize}\qed
	\end{proposition}
	\begin{example}
		Take $f=1132532\in \SF(7)$. Then $\sigma=\chi(f)$ can be obtained as follows:
		\begin{center}   
			$\begin{array}{l}
				(1)\\
				(1,2)\\
				(1,2)(3)\\
				(1,2,4)(3)\\
				(1,2,4)(3)(5)\\
				(1,2,4)(3,6)(5)\\
				(1,2,7,4)(3,6)(5)=\sigma.
			\end{array}$
		\end{center} 
	\end{example}
	The following lemma can easily be deduced from the above proposition and the definition of $\chi$.
	\begin{lemma}\label{lemmaTFE}
		Let $f\in\SF(n)$ and $\sigma=\chi(f)=C_1C_2\cdots C_\ell$. If $\emptyset\neq S\subseteq[n]$, then the following statements are equivalent.
		\begin{enumerate}
			\item $f$ has the property 
			$$\begin{cases}
				f_i=\min(S), \text{ if } i\in S,\\
				f_i\notin S, \text{ else}.
			\end{cases}$$
			\item The elements of $S$ form some cycle $C_i$ in $\sigma$, and the cycle can be written with its elements forming a decreasing sequence. 
			\item $S$ is the underlying set of some cycle $C_i$ with just one weak excedance.\qed
		\end{enumerate}
	\end{lemma}
	
	Consider the bijection $\tau :\SP(n)\mapsto\RGF(n)$ given by $\tau(P)=f$, where $f$ is the canonical form of $P$.
	\begin{definition}
		A Bell permutation of the second kind over $[n]$ is a permutation $\sigma$ obtained from $f\in\RGF(n)$ by applying $\chi$ to $f$, i.\ e., $\sigma=\chi(f)$.
	\end{definition} 
	Let $\BP_2(n):=\chi(\RGF(n))$, the set of all Bell permutations of the second kind over $[n]$, and $\BP_2(n, k):=\{\sigma\in\BP_2(n): \wexc(\sigma)=k\}$.
	
	The restriction of $\chi$ to $\RGF(n)$ is a bijection between $\RGF(n)$ and $\BP_2(n)$. Therefore, $\BP_2(n)$ is a Bell number enumerated set, i.\ e., $|\BP_2(n)|=B(n)$, the $n$'th Bell number.
	
	Since the composition of bijections is a bijection, the map $\lambda=\chi\circ\tau$ is a bijection between $\SP(n)$ and $\BP_2(n)$. 
	\begin{proposition} \label{thmBellchar}
		Let $P=B_1|B_2|\cdots|B_k$ be a set partition, $\sigma$ the permutation $\lambda(P)$, and $C_1C_2\cdots C_\ell$ the cycle decomposition in $\sigma$. Then 
		\begin{enumerate}
			\item $\sigma$ has $k$ weak excedances, 
			\item the set of the weak excedances of $\sigma$ is exactly the interval $[k]=\{1, 2, \ldots, k\}$, and
			\item the set of the minimal elements of the cycles of $\sigma$ is exactly the interval $[\ell]$.
		\end{enumerate}
	\end{proposition}
	\begin{proof}
		The first two items directly follow from Remark \ref{propertyofphitilde1} and the fact that the number of blocks of $P$ is equal to the cardinality of the image set of its canonical form. 
		
		Item 3. By Remark \ref{propertyofphitilde2}, any integer $i\in[n]$ is fixed in $f$ if and only if $i=\min(C_j)$ for some $j$. We show that if $p$ is the maximum fixed point in $f$, then any $q<p$ is also fixed. Suppose that there exist a non-fixed point smaller than $p$. Let $t$ be the maximal of such non-fixed points, i.\ e., the elements of the interval $[t+1,p]$ are all fixed. So $f_t<t$ and $t\notin\{f_1, f_2, \ldots, f_{t+1}=t+1\}$. This implies that $f\notin\RGF(n)$ and this is a contradiction. Therefore, the set of fixed points of $f$ is $[p]$ and hence, $p=\ell$.
	\end{proof}
	The above proposition implies that the distribution of the number of weak excedances on $\BP_2(n)$ is the same as the distribution of the number of blocks on $\SP(n)$, and also that the statistic of the number of cycles on $\BP_2(n)$ has the same distribution as the number of fixed points on $\RGF(n)$. Thus, we have the following.
	\begin{corollary}\label{corcycfx}
		\begin{enumerate}
			\item $|\BP_2(n, k)|=S(n,k)$,
			\item the number of set partitions having $\ell$ blocks containing their own index element is the same as the number of Bell permutations of the second kind having $\ell$ cycles. 
		\end{enumerate}
	\end{corollary}
	The following proposition gives a recursive procedure to check if a permutation is a Bell Permutation of the second kind. We consider the following lemma.
	\begin{lemma}\label{lemmaconcatenate}
		Let $f'\in\SF(n{-}1)$, and let $f\in\SF(n)$ be obtained by concatenating some $j\in[n]$ at the end of $f'$.
		Let $\sigma'=\chi(f')$ and $\sigma=\chi(f)$. If $j\neq n$, then $\sigma$ is obtained from $\sigma'$ by replacing the integer $\sigma'(j)$ by $n$ in $\sigma'$ and appending $\sigma'(j)$ at the end. If $j=n$, then $\sigma$ is obtained by simply appending $n$ at the end of $\sigma'$.\qed
	\end{lemma} 
	\begin{lemma}\label{lemmarecursive}
		A permutation $\sigma=\sigma(1)\sigma(2)\cdots\sigma(n)\in \mathfrak{S}_n$ whose set of weak excedances is an integer interval $[k]$ is in $\BP_2(n)$ if and only if the permutation $\sigma'\in\mathfrak{S}_{n-1}$ obtained from $\sigma$ by replacing the integer $n$ by $\sigma(n)$ in $\sigma|_{[n-1]}$ is in $\BP_2(n-1)$.
	\end{lemma}
	\begin{proof}
		According to Lemma \ref{lemmaconcatenate}, for all permutations $\sigma$, if $f=f_1\cdots f_n = (\chi)^{-1}(\sigma)$ and $\sigma'\in\mathfrak{S}_{n-1}$ is the permutation obtained from $\sigma$ by replacing the integer $n$ by $\sigma(n)$, then the transposition array associated with $\sigma'$ is $f'=f_1f_2\cdots f_{n-1}$. Under the hypothesis that $\wExc(\sigma) = \Image(f)$ is an integer interval $[k]$, the following two conditions are trivially equivalent: 
		\begin{enumerate}
			\item  for all $i \in  [n]$, the set  $\{f_1, f_2, \ldots, f_i\}$ is an integer interval with minimum value $1$.
			\item for all $i \in  [n-1]$, the set  $\{f_1, f_2, \ldots, f_i\}$ is an integer interval with minimum value $1$.
		\end{enumerate}
		That is, $\sigma$ is Bell if and only if $\sigma'$ is Bell. 
	\end{proof}
	For instance, let $\sigma=7245613$. We have $\wExc(\sigma)=[5]$, so $\sigma$ may be a Bell permutation of the second kind. We apply Lemma \ref{lemmarecursive}:
	$\underline{72456}13\rightarrow \underline{32456}1\rightarrow \underline{3245}1\rightarrow \underline{324}1\rightarrow\underline{32}1$. 
	Since $321\in\BP_2(3)$ we can conclude that $\sigma$ and those permutations obtained in the process are Bell permutations of the second kind. But $32541\notin\BP_2(5)$, because $\underline{3254}1\rightarrow \underline{32}1\underline{4}\rightarrow\underline{32}1\in\BP_2(3)$ and $3214\notin\BP_2(4)$.
	
	We give a new proof of the fact that $|\BP_2(n,k)|=S(n,k)$ by showing that the numbers $|\BP_2(n,k)|$ satisfy the recurrence relation of the Stirling number of the second kind. 
	\begin{proposition}\label{thmbellpart}
		The number $|\BP_2(n,k)|$ satisfies the recurrence relation for all positive integers $n, k, n\ge1, 1\le k\le n$:
		\begin{equation}\label{rrstirling2}
			|\BP_2(n,k)| = k|\BP_2(n-1,k)| + |\BP_2(n-1,k-1)|,~ |\BP_2(0,0)|=1.
		\end{equation} 
	\end{proposition}
	\begin{proof}
		We use Lemma \ref{lemmarecursive} to prove the assertion. Any Bell permutation of the second kind $\sigma\in\BP_2(n, k)$ can uniquely be obtained either from a permutation $\sigma'\in \BP_2(n-1, k)$ and an integer $i\in[k]$, or from a permutation $\sigma'\in \BP_2(n-1, k-1)$.  More precisely: if $\sigma'\in \BP_2(n-1, k)$ and $i\in[k]$, then $\sigma$ is obtained from $\sigma'$ by replacing $\sigma'(i)$ by $n$ and then appending $\sigma'(i)$ at the end, i.\ e., $\sigma=\sigma'(i,n)$. In this case $\sigma\in \BP_2(n, k)$, and there are $|\BP_2(n-1,k)|$ possible choices for $\sigma'$ and $k$ possible choices for $i$. Hence this contributes $k|\BP_2(n-1,k)|$ to $|\BP_2(n,k)|$. If $\sigma'\in \BP_2(n-1, k-1)$, then $\sigma$ is obtained from $\sigma'$ by replacing $\sigma'(k)$ by $n$ and then appending $\sigma'(k)$ at the end, i.\ e., $\sigma=\sigma'(k,n)$. In this case $\sigma\in \BP_2(n, k)$, and $\sigma$ has $|\BP_2(n-1,k-1)|$ possibilities. By combining the two cases we have (\ref{rrstirling2}).
	\end{proof}
	Let $P\in\SP(n,k)$ and $\MX(P)=\{\max(B_i): 1\le i\le k\}$. By the above proposition, Remark \ref{propertyofphitilde1}, and the fact that the maximum elements of the blocks of $P$ are the rightmost occurrences in $\tau(P)$ we have the following corollaries.
	\begin{corollary}
		We have $\MX(P)=\wExcL(\sigma)$, where $\sigma=\lambda(P)$.
	\end{corollary}
	\begin{corollary}
		The bistatistics $(\blocks,\fx)$ on the set $\SP(n)$ has the same distribution as $(\wexc,\cyc)$ on the set $\BP_2(n)$.
	\end{corollary}
	\begin{remark} 
		The cardinality of the set $\BP_2(n,n{-}1)$ is equal to the number $S(n,n{-}1)$ of set partitions over $[n]$ having $n{-}1$ blocks, which, as is well known, is  equal to $\binom{n}{2}$.
	\end{remark}
	OEIS entry number \seqnum{A259691} presents the sequence of the numbers $T(n{-}1,\ell)$, counting set partitions over $[n]$ where exactly $\ell$ blocks contain their own index element. These numbers satisfy the relation:
	\begin{equation}\label{eqncycle}
		T(n{-}1,\ell)=\sum_{i=0}^{n-\ell}{n-\ell\choose i}\ell^{i+1}B(n-\ell-i),
	\end{equation}
	where $T(n{-}1,n)=1$. Thus, by Corollary \ref{corcycfx} the number of Bell permutations of the second kind over $[n]$ having exactly $\ell$ cycles is also equal to $T(n{-}1,\ell)$. 
	
	We can refine (\ref{eqncycle}) by adding an additional parameter counting the number of weak excedances of the permutation.
	Consider a permutation $\sigma'\in\BP_2(n{-}1,k)$, an integer $i\in[k+1]$ representing a weak excedance and the permutation $\sigma=\sigma'(i,n)$, the product of $\sigma'$ and the transposition $(i,n)$. Then, by Proposition \ref{propcycsigma}, the numbers of cycles of $\sigma$ and $\sigma'$ are equal, except when $i=k+1=n$ (i.\ e., both $\sigma'$ and $\sigma$ are the identity permutations) and $i=k+1$, in which case $\cyc(\sigma)=\cyc(\sigma')+1$. Thus, the number $T(n{-}1,k,\ell)$ of Bell permutations of the second kind over $[n]$ having exactly $k$ weak excedances and $\ell$ cycles satisfies:
	\begin{equation}\label{bellwexccyc}
		T(n{-}1,k,\ell)=\begin{cases} \delta_{k,n}, \text{ if } \ell=n,\\
			\sum\limits_{i=0}^{n-\ell}{n-\ell\choose i}\ell^{i+1}S(n{-}\ell{-}i,k{-}\ell), \text{ else. }\end{cases}
	\end{equation}
	where $\delta_{*,*}$ is the Kronecker delta function.
	Therefore, we have the following.
	\begin{proposition} For $n\ge 1$ we have
		\begin{equation} 
			\sum\limits_{\sigma{\in}\BP_2(n)}x^{\wexc(\sigma)}y^{\cyc(\sigma)}=\sum\limits_{k=1}^{n}\sum\limits_{\ell=1}^kT(n{-}1,k,\ell)x^ky^\ell.
		\end{equation}\qed
	\end{proposition}
	\begin{corollary}
		The number of Bell permutations of the second kind over $[n]$ having exactly $1$ cycle equals $B(n{-}1), n\geq1$. 
	\end{corollary}
	\begin{proof}
		Let $\ell=1$ in (\ref{bellwexccyc}) and take the sum over all $1\le k\le n$.
	\end{proof}
	\subsection{A bijection between Bell permutations of the first and the second kind}
	\label{bijBellperm1and2}
	In this subsection, we present a bijection between the set $\BP_1(n)$ of Bell permutations introduced by Poneti and Vajnovszki \cite{Po-Va} (which we will call \emph{Bell permutations of the first kind}) and the set $\BP_2(n)$ of Bell permutations of the second kind.
	
	First, we recall the definition of Bell permutations of the first kind. Let $P=B_1|B_2|\cdots |B_k$ be a set partition over $[n]$ in its standard representation and let $\mu : \SP(n)\mapsto \BP_1(n)$, where the permutation $\mu(P)$ is constructed as follows:
	\begin{itemize}
		\item reorder all integers in each block $B_i$ in decreasing order;
		\item transform each of these blocks into a cycle.
	\end{itemize}
	For instance, if $P=1279|356|48$, then $\mu(P)=(9,7,2,1)(6,5,3)(8,4)$.
	
	By Lemma \ref{lemmaTFE}, if $\mu(P)=\sigma\in \BP_1(n)$ and $f=\chi^{-1}(\sigma)$ is its transposition array, then for all $i \in [n]$, $$f_i=\mbox{ minimum of the block of $P$ containing } i.$$
	Recall also that if $\sigma\in \BP_2(n)$ and $f=\chi^{-1}(\sigma)=\tau(P)$ is its transposition array, then for all $i \in [n]$, $$f_i=\mbox{ index of the block of } P \mbox{ containing } i.$$
	
	Thus, we have the bijection $\beta :=\lambda\circ\mu^{-1}: \BP_1(n)\mapsto \BP_2(n)$. As we shall see that $\beta$ can be described concretely as follows.
	\begin{proposition} \label{thmbeta}
		Let $\sigma=C_1C_2\cdots C_k\in\BP_1(n)$, written in cycle notation, where each cycle is ordered decreasingly. Let $\sigma'$ be constructed from $\sigma$ according to the rule:
		for $i~=~k, k{-}1, \ldots, 2$, if the integer $i$ is not in the $i$'th cycle, then insert the sequence of elements of the $i$'th cycle after $i$ in the cycle containing $i$. Then $\sigma'=\beta(\sigma)$.
	\end{proposition}
	\begin{proof}
		Let $f=\chi^{-1}(\sigma)$, $\nu$ be the transformation that normalizes $f$ via the order-preserving bijection of $\Image(f)$ into $[\IMA(f)]$, and $f'=\nu(f)$. Let $\gamma=\chi\circ\mu^{-1}$. Observe that $\nu=\gamma\circ\tau$. By Lemma \ref{lemmaInom} and inspection we have 	$\chi^{-1}(\sigma')=f'=\nu(f)=\tau\circ\mu^{-1}(\sigma)$. Indeed, if $f^{(i)}$ is the transposition array associated with the permutation obtained after the $i$'th step of the procedure, then it can be verified that for all integers $j$ in the cycle $C_i$ one has $f^{(i)}(j)=i$ and the image of such integers $j$ does not change in the following steps. In other words, the following diagram is commutative.
		\begin{figure}[H]
			\centering
			\includegraphics[scale=.80]{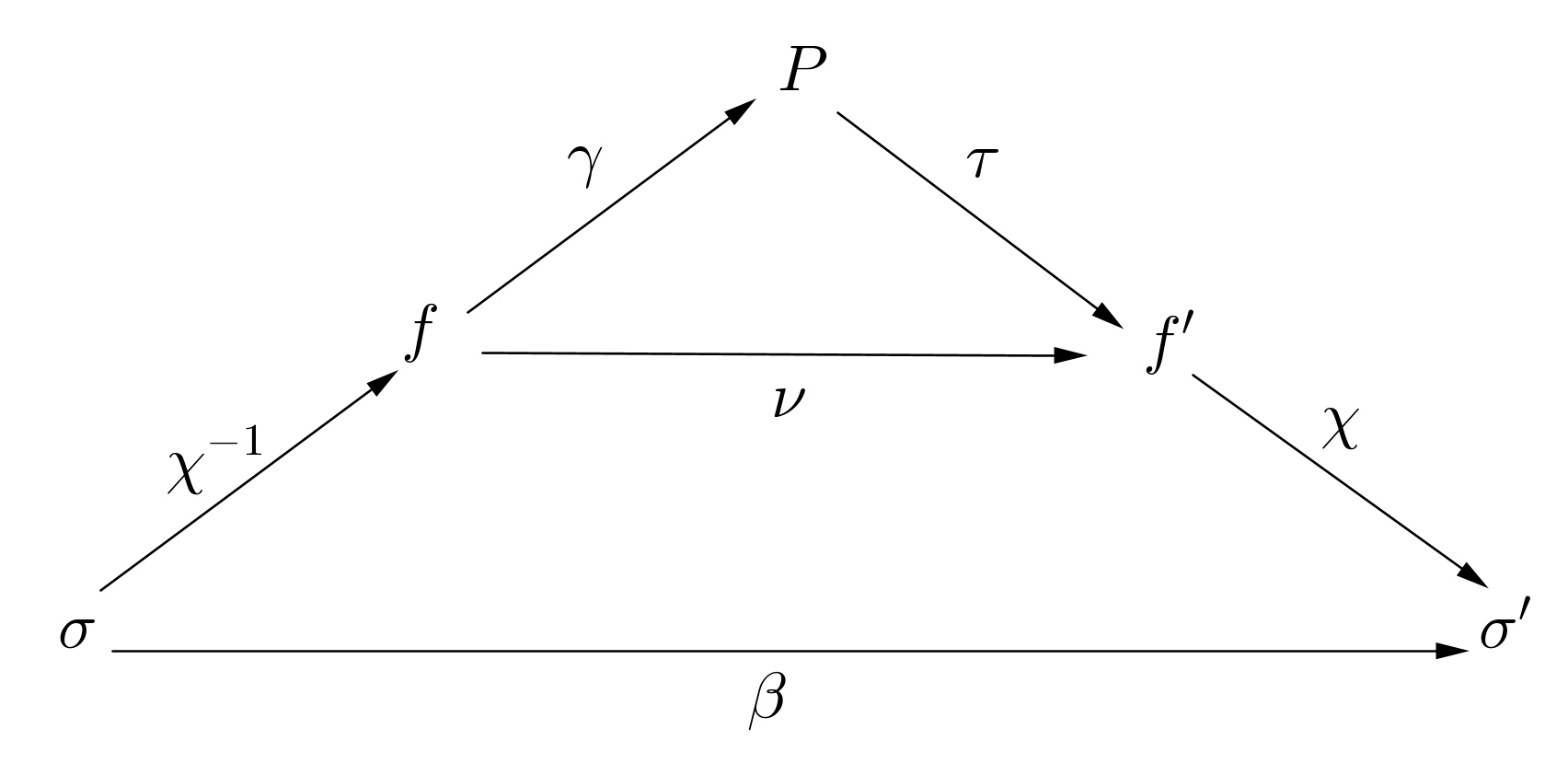}
			\caption{}
		\end{figure}
		So we have $ \sigma'=\chi(f')=\chi \circ \tau \circ \mu^{-1}(\sigma)=\lambda\circ\mu^{-1}(\sigma)=\beta(\sigma)$ indeed. 
	\end{proof}
	For instance, let $\sigma=(9,7,2,1)(6,5,3)(8,4)$. Then $\sigma'$ is obtained as:
	$$\sigma=(9,7,2,1)(6,5,3)(8,4)\longrightarrow(9,7,2,1)(6,5,3,8,4)\longrightarrow(9,7,2,6,5,3,8,4,1)=\sigma'.$$		

	We can also describe directly $\vartheta:=\beta^{-1}$ as follows. Take  $\sigma'\in \BP_2(n)$  and let $C_1C_2\cdots C_l$ be its cycle decomposition. Assume that $\wExc(\sigma)=[k]$. For $i= 2, \ldots, k$, if $i$ is not the minimum of its own cycle $C_j$, then form a new cycle by taking out of $C_j$ the longest sequence of integers greater than $i$ starting immediately after $i$, and modify the cycles. The resulting permutation is $\sigma=\vartheta(\sigma')$.
	For instance, let $\sigma=468912357=(\mathbf{1,4,}9,7,\mathbf{3,5,}8)(\mathbf{2,}6)$ in cycle notation and with the weak excedances in bold. Then $\sigma$ is obtained as:
	$$\sigma'=(1,4,9,7,3,5,8)(\mathbf{2,}6)\longrightarrow(1,4,9,7,\mathbf{3})(2,6)(5,8)\longrightarrow(1,\mathbf{4,}3)(2,6)(5,8)(9,7)\longrightarrow$$$$(1,4,3)(2,6)(\mathbf{5})(9,7)(8)=\sigma~~~~~~~~~~~~~~~~~~~~~~~~~~~~~~~~~~~~~~~~~~~~~~~~~~~~~~~~~~~~$$
\begin{remark}
	Under the bijection $\beta: \sigma\mapsto\sigma'$, the number of cycles of $\sigma$ is equal to the number of weak excedances of $\sigma'$.
\end{remark}

\medskip
The OEIS entry number \seqnum{A026898} enumerates the number of set partitions over $[n+1]$ whose minima form an interval of positive integers starting with $1$. By Corollary \ref{corcycfx} and Proposition \ref{thmbellpart}, these set partitions correspond to Bell permutations of the second kind over $[n+1]$ having equal number of weak excedances and number of cycles. Also notice that $\BP_1(n)\cap\BP_2(n)=\{\sigma\in\BP_2(n): \wexc(\sigma)=\cyc(\sigma)=\ell\}$. Thus and by (\ref{bellwexccyc}) we have the following
\begin{corollary}
	For $n\ge 1$
	\begin{equation}
		|\BP_1(n)\cap\BP_2(n)|=1+\sum\limits_{\ell=1}^n\ell^{n-\ell+1}. 
	\end{equation}\qed
\end{corollary}
\section{Other classes of Bell enumerated subexcedant functions}
\label{otherclassSef}
In this section we present two families of subexcedant functions also counted by the Bell numbers. 

Let $f=f_1f_2\cdots f_n\in\SF(n)$. Recall that $i$ is a leftmost occurrence in $f$ if $f_i~\notin~\{f_1, \ldots, f_{i-1}\}$, where $i\in[n]$. Clearly $1$ is a leftmost occurrence. We say that $i>1$ is a \textit{repetition} in $f$ if it is not a leftmost occurrence. 

A subexcedant function $f$ is said to avoid a pattern $212$ (or $121$) if there do not exist some indices $a<b<c$ such that $f_a=f_c>f_b$ (or $f_a=f_c<f_b$, respectively).

The first family we consider is the set $\SF_1(n)$ of subexcedant functions over $[n]$ such that for $j\in\Image(f)$, the set of all $f^{-1}(j)$ form an integer interval. For instance, $1133222\in\SF_1(7)$. The following remark characterizes the set $\SF_1(n)$ in terms of pattern avoidance.
\begin{remark}
	A subexcedant function $f\in\SF_1(n)$ if and only if $f$ is $212$ and $121$-avoiding.\qed
\end{remark}
We let $\SF_1(n,k):=\{f\in\SF_1(n): \IMA(f)=k\}$. Define the map $\omega:~\SF_1(n,k)\mapsto\SP(n,n+1-k)$ by $\omega(f)=P$, where $P$ is the set partition obtained from $f$ as follows: initialize the first block with $f_1=1$ as a minimum, the remaining $n-k$ blocks with the repetitions as minima, and finally insert a leftmost occurrence $i>1$ in the $j$'th block, where $j=|[f_i]\backslash\{f_1, f_2, \ldots, f_{i-1}\}|$.
\begin{example}
	Consider $f=111334268\in\SF_1(9,6)$. The set of repetitions of $f$ is $\{2, 3, 5\}$. So there are $4$ blocks initialized as: $1\cdots|2\cdots|3\cdots|5\cdots$. Since $|[3]\backslash\{1\}|=2$, the leftmost occurrence $4$ is inserted in the $2$-nd block. From $|[4]\backslash\{1,3\}|=2$ we determine that $6$ is inserted in the $2$-nd block, and so on. Thus, we obtain the set partition $\omega(f)=P=1~7|2~4~6~8|3~9|5$. Observe that $P\in\SP(9,4)$.
\end{example}
Conversely, assume that the values $f_1{=}1, f_2, \ldots, f_{i-1}$ have already been computed. If $i$ is in the $j$'th block of $P$ and $i>\min(B_j)$, then let $f_i$ be the $j$'th smallest element of the set $[n]\backslash\{f_1, f_2, \ldots, f_{i-1}\}$. If $i=\min(B_j), j>1$, then let $f_i=f_{i-1}$. It is easy to see that $f=\omega^{-1}(P)$. 
\begin{proposition}
	The map $\omega$ is a bijection.\qed
\end{proposition}	
\begin{corollary}
	For $n\geq1$, we have $$|\SF_1(n,k)|=S(n,n+1-k),$$
	where $S(n,k)$ is the Stirling number of the second kind.\qed
\end{corollary}

The second family of subexcedant functions we consider is given as follows. 

For $f\in \SF(n)$, we define $\RmL(f)$ to be the \emph{subword} of the rightmost letters of $f$ in the order they appear in $f$, i.\ e, if $f=f_1f_2\cdots f_n$, then $\RmL(f)$ is the subword of $f$ composed of all $f_i$'s such that $i$ is a rightmost occurrence of $f$. Note that $\RmL(f)=\Image(f)$ as sets. For instance, if $f=121135623$, then $\RmL(f)=15623$. Recall that the rightmost letters of $f$ correspond to the weak excedances of the corresponding permutation $\sigma=\chi(f)$. Thus, the subword $\RmL(f)$ is increasing if and only if the subword of weak excedance letters of $\sigma$ is increasing. For the function $f=121135623$, the corresponding permutation is $\sigma=489367125$. The subword of its weak excedance letters is $48967$ and it is not increasing. 

We let $\SF_2(n)$ denote the set of subexcedant functions over $[n]$ whose subword of the rightmost letters is increasing. Also, let $\SF_2(n,k):=\{f\in\SF_2(n) : \IMA(f)=k\}$.
\begin{theorem}
	The number of permutations in $\mathfrak{S}_n$ having increasing subword of weak excedance letters is the $n$'th Bell number $B(n)$.
\end{theorem}
\begin{proof}
	We give two proofs via the transposition arrays of such permutations. We first prove directly that the cardinality of the set $\SF_2(n,k)$ is equal to $S(n,k)$, which satisfies the relation in (\ref{rrstirling2}), and then provide another proof by presenting a bijection between $\SF_2(n)$ and $\RGF(n)$.
	
	Suppose that the subword of weak excedance letters of a permutation $\sigma$ is increasing. Let $f$ be the transposition array of $\sigma$, i.\ e., $f=\chi^{-1}(\sigma)$ with $\RmL(f)=f_{i_1}f_{i_2}\ldots f_{i_k}$. Then, we have $f_{i_1}<f_{i_2}<\cdots <f_{i_k}$ and $i_1<i_2<\cdots<i_k$. Therefore, $f\in\SF_2(n,k)$. Each such subexcedant function can be obtained in either of the following ways. Consider a subexcedant function $f\in\SF_2(n{-}1,k)$. Let $a$ be an element of $\Image(f)$, and let $f'$ be obtained from $f$ by inserting the value $a$ in the position $a$. Then $f'\in\SF_2(n,k)$ and $\RmL(f')=\RmL(f)$. Since there are $k$ possible choices for $a$, this contributes $k|\SF_2(n-1,k)|$ to the number $|\SF_2(n,k)|$. Consider a subexcedant function $f\in\SF_2(n{-}1,k{-}1)$ with $\RmL(f)= f_{i_1}<f_{i_2}< \cdots <f_{i_{k-1}}$, where $i_1<i_2<\cdots<i_{k-1}$. Let $f'$ be obtained from $f$ by appending $n$ at its end. Then $f'\in\SF_2(n,k)$ and $\RmL(f')=\langle f_{i_1}<f_{i_2}< \cdots <f_{i_{k-1}}<n\rangle$, where $i_1<i_2<\cdots<i_{k-1}<n$. This contributes $|\SF_2(n{-}1,k{-}1)|$ to the number $|\SF_2(n,k)|$. Hence, by combining the cases we have the proof.
	
	Alternatively, we present a bijection between the sets $\SF_2(n)$ and $\RGF(n)$. Let $f\in\SF_2(n)$ and $f'$ be the function obtained from $f$ as follows. For $i=n, n-1, \ldots, 2, 1$: let $g^{(n)}=f$ and $g^{(i)}$ be the function obtained from $g^{(i+1)}$ by deleting the largest fixed point. Note that $g^{(i)}$ is a subexcedant function over $[i]$. Now let $j$ be the largest fixed point in the function $g^{(i)}$, set $f_i'=j'$, where $j'$ is the normalized value of $j$ under the map $\nu$ given in Proposition \ref{thmbeta}. We note that $f'$ is a restricted growth function, and that $\IMA(f)=\IMA(f')$. 
	
	Conversely, let $f'\in\RGF(n)$. We obtain $f$ uniquely from $f'$ as follows. Suppose that the function $g^{(i-1)}$ has already been computed. This is a subexcedant function over $[i{-}1]$. Then at the $i$'th step: if $j=f_i'\le\IMA(g^{(i-1)})$, and $a$ is the $j$'th smallest element in $\Image(g^{(i-1)})$, then insert $a$ also as a value in the function $g^{(i-1)}$ in the position $a$; otherwise, let $g^{(i)}_i=i$. It can easily be seen that $f=g^{(n)}\in\SF_2(n)$. Therefore, $f\mapsto f'$ is a bijection.
\end{proof}
\begin{example}
	Take ${f=11131338\in\SF_2(8)}$. Then ${\IMA(f)=3}$ and the corresponding $\RGF$ $f'=f_1'f_2'\cdots f_8'$ is obtained as follows.
	$$g^{(8)}=11131338, ~f_8'=3$$
	$$g^{(7)}=1113133, ~f_7'=1$$
	$$g^{(6)}=113133, ~f_6'=2$$
	$$g^{(5)}=11133, ~f_5'=1$$
	$$g^{(4)}=1133, ~f_4'=2$$
	$$g^{(3)}=113, ~f_3'=2$$
	$$g^{(2)}=11, ~f_2'=1$$
	$$g^{(1)}=1, ~f_1'=1$$
	Therefore, $f'=11221213\in\RGF(8)$.
\end{example}

\section{Bijections on set partitions}
\label{involution}
In this section, we present some bijections on set partitions. In particular, we present an involution that interchanges the number of merging blocks (that we define below) and the number of successions. We use these bijections to study the power series generating function for the distribution of these statistics, and to deduce some structural results for set partitions. 

For $ n\ge1$, we shall describe a partition of $\SP(n)$ into equivalence classes. Set partitions within each class are closely related. Each class will contain exactly one merging-free partition. Since there are exactly $B(n{-}1)$ merging-free partitions, the same is true for the number of classes. The size of each class is a power of two.

Recall that a set partition $P=B_1|B_2|\cdots|B_k\in\SP(n,k)$ in standard form satisfies the condition $\min(B_i)<\min(B_{i+1}), 1\leq i<k$. 

\subsection{Merging and successions equivalence}
In this subsection, we discuss how to transform a merging block of a set partition into a succession and vice versa. 

Let $\mathcal{T}_n^a:=\{P\in\SP(n): a\in\MB(P)\}$ and $\mathcal{R}_n^a:=\{P\in\SP(n): a\in\Succeq(P)\}$. We always assume that $a\in[2,n]$. Further for any $A\subseteq[2,n]$, let $\mathcal{T}_n^A:=\{P\in\SP(n): \MB(P)=A\}$ and $\mathcal{R}_n^A:=\{P\in\SP(n): \Succeq(P)=A\}$. It can easily be seen that $\mathcal{T}_n^a=\bigcup\limits_{\substack{A\\a\in A}}\mathcal{T}_n^A$, similarly  $\mathcal{R}_n^a=\bigcup\limits_{\substack{A\\a\in A}}\mathcal{R}_n^A$.

\begin{remark}We recall (\cite[Proposition 1.1]{Be-Ma}) \label{remmerfree}
	that the number $|\mathcal{T}_n^\emptyset|$ of merging-free partitions over $[n]$ equals the  Bell number $B(n{-}1), n\ge1$. Likewise, the sequence of the number of set partitions over $[n]$ having $m$ successions is presented in OEIS entry number \seqnum{A056857} (also, see Munagi \cite{Mu}).
\end{remark}

We define the operation $\Swap_a^{(i,j)}$ on a set partition $P=B_1|B_2|\cdots |B_k$, where $i$ and $j$ are two integers in $[k]$ and $a\in [n]$. If $i=j$ or $a\notin B_i\cup B_j$, we let $\Swap_a^{(i,j)}(P)=P$. Else, we let $I_a$ be the maximal integer interval in $B_i\cup B_j$ that starts with $a$, and we move the elements of $I_a$ lying in $B_i$ to $B_j$ and vice versa. 

For instance, let $P=1~3~4~6~8|2~5~9|7$. Then $\Swap_3^{(1,2)}(P)=1~5~8|2~3~4~6~9|7$,  $\Swap_7^{(1,3)}(P)=1~3~4~6~7|2~5~9|8$, $\Swap_3^{(1,1)}(P)=\Swap_3^{(2,3)}(P)=P$, and $\Swap_7^{(2,3)}(P)=1~3~4~6~8|2~5~7~9| ~$, with the last new block empty. (Here, strictly speaking, $\Swap_7^{(2,3)}(P)$ is not a set partition. However, in our applications of $\Swap$, such an empty block never appears.)

We now define the following maps.
\begin{enumerate} 
	\item Consider $P=B_1|B_2|\cdots |B_k\in\mathcal{T}_n^a$. Then $a{=}\min(B_i)$ and $a{-}1\in B_j$ for certain $i$ and $j$. Note that then $\min(B_j)\le a{-}1<\min(B_i)$, whence $j<i$. Define the map $\mu_a: \mathcal{T}_n^a\mapsto \mathcal{R}_n^a$ by $\mu_a(P)=P'$, where $P'$ is obtained from $P$ as follows. Let $P^*$ be the set partition obtained by merging the blocks $B_{i-1}$ and $B_i$, and put $P'=\Swap_a^{(i-1,j)}(P^*)$.  
	We note that $a$ becomes a succession of $\mu_a(P)$. Later we will show that $\NMB(P')=\NMB(P)$. For instance, let $P{=}1~3~5~7~10|2~4|6~8|9$. We have $\MB(P)=\{6, 9\}, \Succeq(P)=\emptyset$. If $a=6$, then $i=3, j=1$ and $P^*=1~3~5~7~10|2~4~6~8|9$. Thus, $P'=\mu_6(P)=1~3~5~6~8~10|2~4~7|9\in\mathcal{R}_{10}^6$. Note that $\MB(P')=\{9\}, \Succeq(P')=\{6\}$, and that $\NMB(P')=\NMB(P)=\{1, 2\}$.
	\item Consider $P=B_1|B_2|\cdots |B_k\in\mathcal{R}_n^a$. Then $a-1, a\in B_i$ for some $i$. Define the map $\rho_a: \mathcal{R}_n^a\mapsto \mathcal{T}_n^a$ by $\rho_a(P)=P'$, where $P'$ is obtained from $P$ as follows. Let $j$ be the smallest positive integer such that the elements $1, 2, \ldots, a-1$ are in the first $j$ blocks of $P$. Apply $\Swap_a^{(i,j)}$ to $P$, and then split the modified block $B_j$ before $a$. We note that the succession $a$ becomes the minimum element of a merging block of $\rho_a(P)$. Further, $\NMB(P')=\NMB(P)$. For instance, let $P=1~3~4~6~9|2~5~8|7|10$ with $\MB(P)=\{10\}, \Succeq(P)=\{4\}$, and let $a=4$. So $i=1, j=2$, and $\Swap_4^{(1,2)}(P)=1~3~5~9|2~4~6~8|7|10$. Hence, we have $\rho_4(P)=P'=1~3~5~9|2|4~6~8|7|10\in\mathcal{T}_8^4$. Observe that $\MB(P')=\{4, 10\}, \Succeq(P')=\emptyset$, and that $\NMB(P')=\NMB(P)=\{1, 2, 7\}$.
\end{enumerate}
\begin{lemma}\label{lemmaMbtoSuc}
	\begin{enumerate}
		\item If $a\in\MB(P)$ and $P'=\mu_a(P)$, then $\MB(P')=\MB(P)\backslash\{a\}$,  $\Succeq(P')=\Succeq(P)\cup\{a\}$, and $\NMB(P')=\NMB(P)$.
		\item If $a\in\Succeq(P)$ and $P'=\rho_a(P)$, then $\Succeq(P')=\Succeq(P)\backslash\{a\}$ and $\MB(P')=\MB(P)\cup\{a\}$, and $\NMB(P')=\NMB(P)$.
	\end{enumerate}
\end{lemma}
\begin{proof}
	We provide only the proof of the former item since the proof of the latter would be analogous. 
	
	Let $P=B_1|\cdots|B_k\in\mathcal{T}_n^a$, where $a \in B_i,~ a{-}1\in B_j$ for some $j<i\le k$. So $\max(B_{i-1})<\min(B_i)$ since $B_i$ is merging. Let $I_a$ denote the interval of integers moved by $\mu_a$ (by this we mean the interval moved by the $\Swap$ operation in the procedure of $\mu_a$). Let $P'=\mu_a(P)=B_1'|\cdots|B_{k-1}'$. We consider two cases.
	
	If $j=i{-}1$, then $B_x'=B_x$ for $x<i{-}1$, $B_{i-1}'=B_{i-1}\cup B_i$, and $B_x'=B_{x+1}$ for $i\le x<k$. This implies $\max(B_{i-1}')=\max(B_i)$ and $\min(B_{i-1}')=\min(B_{i-1})$. So $B_{i-1}'$ (resp. $B_i'$) is merging if and only if $B_{i-1}$ (resp. $B_{i+1}$) is merging. Thus, $\MB(P')=\MB(P)\backslash\{a\}, \Succeq(P')=\Succeq(P)\cup\{a\}$. 
	
	Now suppose that $j<i{-}1$. In this case $B_x'=B_x$ for $x<j$ or $j<x<i{-}1$ and $B_x'=B_{x+1}$ for $i\le x<k$, and $\max(B_j)>\min(B_{j+1})$ since $\max(B_j)\ge a{-}1$ and $\min(B_{j+1})<a$. Since the integers of the interval $I_a$ are greater than or equal to $a$ and $\mu_a$ swaps these integers between $B_j$ and $B_i$, $\max(B_j')\ge \max(B_j)$. Observe that $\min(B_{j+1}')=\min(B_{j+1})$. Thus, $\max(B_j')>\min(B_{j+1}')$. Further, we have that $\max(B_{i-1}')\le\max(B_i)$ and $\min(B_i')=\min(B_{i+1})$. Hence $\max(B_{i-1}')>\min(B_{i}')$ if and only if $\max(B_{i})>\min(B_{i+1})$. Therefore, no new merging block is created in this process and hence $\MB(P')=\MB(P)\backslash\{a\}$. 
	
	On the other hand, let us show that the process does not create any new succession other than $a$. If $b{-}1, b\in I_a, b>a$, then either both of them belong to the same block in $P$ and thus $\mu_a$ moves them together to the other block, or they belong to different blocks and thus $\mu_a$ swaps them. Thus $\Succeq(P')=\Succeq(P)\cup\{a\}$. 
	
	Furthermore, observe that neither $\mu_a$ nor $\rho_a$ moves the minimum element of a non-merging block. Thus $\NMB(P)$ is preserved under these maps.
\end{proof}
\begin{lemma}\label{lemmaInverse}
	We have $\rho_a\circ\mu_a=id_{\mathcal{T}_n^a}$ and $\mu_a\circ\rho_a=id_{\mathcal{R}_n^a}$. In other terms, $\mu_a$ and $\rho_a$ are inverses of each other.
\end{lemma}
\begin{proof} 
	Since a succession cannot be the minimum element of a block for any set partition $P$, we have always $\MB(P)\cap\Succeq(P)=\emptyset$. We first prove that $\rho_a\circ\mu_a=id_{\mathcal{T}_n^a}$. Let $P=B_1|B_2|\cdots |B_k\in\mathcal{T}_n^a$, suppose that $a\in B_{i}$ and $a{-}1\in B_j$ with $j<i\le k$. Since $B_i$ is merging and $P$ is in standard form, $B_{i-1}\subseteq[a{-}1]\subseteq\cup_{\ell=1}^{i-1}B_\ell$. Let $I_a$ be the maximal integer interval moved by $\mu_a$, so $I_a\subseteq B_{i-1}\cup B_i\cup B_j$. After applying $\mu_a$ the integer $a$ becomes a succession in $P'=\mu_a(P)=B_1'|B_2'|\cdots |B_{k-1}'$, i.\ e., $a{-}1, a\in B_j'$, and since $\mu_a$ merges the blocks $B_i$ and $B_{i-1}$, the block $B_{i-1}'$ in $P'$ is the rightmost block containing some integer(s) smaller than $a$. Therefore, when $\rho_a$ is applied to $P'$ it splits precisely this block to create a merging block. Thus, if $I_a'$ is the maximal integer interval moved by $\rho_a$ (i.\ e., by the $\Swap_a^{(i-1,j)}$), then we have $I_a'= I_a$ because $B_{i-1}'\cup B_j'=B_{i-1}\cup B_i\cup B_j$. Therefore, $\rho_a$ reverses the action of $\mu_a$ and $\rho_a\circ\mu_a$ is the identity on $\mathcal{T}_n^a$.
	
	Next we prove that $\mu_a\circ\rho_a=id_{\mathcal{R}_n^a}$. Suppose that $P=B_1|B_2|\cdots |B_k\in\mathcal{R}_n^a$. If $\rho_a$ breaks a succession $a\in B_i$ for some $i\leq k$,	and creates a merging block, say $B_{j}'$ for some $j$, in $P'=\rho_a(P)=B_1'|B_2'|\cdots |B_{k+1}'$, then $a{-}1\in B_j'$ and the interval of integers moved by $\mu_a$ is the same as the interval moved by $\rho_a$. So, the map $\mu_a$ reverses the action of $\rho_a$ and hence $\mu_a\circ\rho_a=id_{\mathcal{R}_n^a}$.
\end{proof}

\begin{lemma}\label{propcommuting}
	For any $a\neq b\in[2,n]$, we have
	\begin{enumerate}
		\item $\mu_a\circ\mu_b=\mu_b\circ\mu_a$ on $\mathcal{T}_n^a\cap\mathcal{T}_n^b$, 
		\item $\rho_a\circ\rho_b=\rho_b\circ\rho_a$ on $\mathcal{R}_n^a\cap\mathcal{R}_n^b$,  and
		\item $\mu_a\circ\rho_b=\rho_b\circ\mu_a$ on $\mathcal{T}_n^a\cap\mathcal{R}_n^b$.
	\end{enumerate}
\end{lemma}
\begin{proof}
	Item 1. Suppose that $P=B_1|B_2|\cdots |B_k\in\mathcal{T}_n^a\cap\mathcal{T}_n^b$ and assume, without loss of generality, that $a=\min(B_{i_1})<b=\min(B_{i_2})$ for some $i_1<i_2\le k$. Let $I_a:=I_{a,P}$ and $I_b:=I_{b,P}$ be the maximal integer intervals moved by $\mu_a$ and $\mu_b$ in $P$, respectively. Let $a{-}1\in B_j$ for some $j<i_1$. Then $I_a\subseteq B_{i_1-1}\cup B_{i_1} \cup B_j$. 
	
	Suppose that $b{-}1\notin I_a$. Then $\mu_a$ does not move $b{-}1$. Observe then that $I_{a,P}$ is a subset of the maximal integer interval $I_{a,\mu_b(P)}$ moved by $\mu_a$ in $\mu_b(P)$. Let $\alpha\notin I_{a,P}$ be the smallest integer greater than $a$. If $\alpha\in I_{a,\mu_b(P)}$, then $\alpha$ would be in the $j$'th block of $\mu_b(P)$. Since $b{-}1\notin I_{a,P}$ (whence $\alpha\le b{-}1$) and $\mu_b$ has only moved  integers greater than $b{-}1$, we would have $\alpha\in B_j$. Therefore, instead, $I_{a,\mu_b(P)}\subseteq I_{a,P}$, and hence $I_{a,\mu_b(P)}=I_{a,P}$. Similarly, $I_{b,\mu_a(P)}=I_{b,P}$. Therefore, we have $\mu_a\circ\mu_b=\mu_b\circ\mu_a$. 
	
	We now suppose that $b{-}1\in I_a$. Then $b{-}1$ is either in the block $B_j$ or in the block $B_{i_1}$ of $P$. In either case $I_a=[a,b{-}1]\subseteq B_{i_1-1}\cup B_{i_1} \cup B_j$, thus we have $i_1{+}1=i_2$, i.\ e., the block containing $a$ and the block containing $b$ in $P$ are adjacent. First, in addition, assume that $b{-}1\in B_j$:
	$$P=B_1|\cdots|\underbrace{\cdots a{-}1\,\cdots b{-}1\,\cdots}_{B_j}|\cdots|B_{i_1-1}|\underbrace{a\cdots}_{B_{i_1}}|\underbrace{b\cdots}_{B_{i_1+1}}|\cdots|B_k.$$
	Consider the product $\mu_b\circ\mu_a$. If $P'=\mu_a(P)=B_1'|B_2'|\cdots |B_{k{-}1}'$, then $b{-}1\in B_{i_1-1}'$ and $b\in B_{i_1}'=B_{i_1+1}$:
	$$P'=B_1'|\cdots|\underbrace{\cdots a{-}1\,a\cdots}_{B_j'}|\cdots|\underbrace{\cdots b{-}1}_{B_{i_1-1}'}|\underbrace{b\cdots}_{B_{i_1}'}|\cdots|B_{k{-}1}'.$$
	Now when $\mu_b$ is applied to $P'$, it simply merges the block $B_{i_1}'$ to $B_{i_1-1}'$ because $I_b\subseteq B_{i_1}'$. Then we have $P''=\mu_b(P')=B_1''|B_2''|\cdots |B_{k-2}''$ as follows. 
	$$P''=B_1''|\cdots|\underbrace{\cdots a{-}1\,a\cdots}_{B_j''}|\cdots|\underbrace{\cdots b{-}1\,b\cdots}_{B_{i_1-1}''}|\cdots|\cdots|B_{k{-}2}''.$$
	Now consider the product $\mu_a\circ\mu_b$. Since $b{-}1\in B_j, b\in B_{i_1+1}$, the interval $I_b\subseteq B_{i_1+1}\cup B_{j}$. Thus $\mu_b$ move the elements of $I_b$ and merges the modified block $B_{i_1+1}$ with the block $B_{i_1}$, i.\ e., we obtain a set partition $P^*=\mu_b(P)$:
	$$P^*=B_1^*|\cdots|\underbrace{\cdots a{-}1\cdots b{-}1\,b\cdots}_{B_j^*}|\cdots|B_{i_1-1}^*|\underbrace{a\cdots}_{B_{i_1}^*}|\cdots|B_{k{-}1}^*.$$
	Since $\mu_a$ moves $b{-}1$, in this case when $\mu_a$ is applied to $P^*=\mu_b(P)$, the interval $I_{a,P^*}=I_{a,P}\cup I_{b,P}$. So $\mu_a$ restores those elements that $\mu_b$ moved from $B_j$ to $B_{i_1+1}$ back to the $j$'th block of $\mu_b(P)$ and vice-versa. Therefore, $\mu_a(\mu_b(P))=\mu_a(P^*)=P''$. 
	
	In the subcase where $b{-}1 \in I_a$ and $b{-}1\in B_{i_1}$ the argument is similar. Hence $\mu_a\circ\mu_b=\mu_b\circ\mu_a$ in all cases. 
	
	For Item 2, and Item 3, we use the equality in Item 1, and the fact that $\mu_a$ and $\rho_a$ are inverses (Lemma \ref{lemmaInverse}). So
	\begin{align*}
		\rho_a\circ\rho_b&=\rho_b\circ\rho_a\circ\mu_a\circ\mu_b\circ\rho_a\circ\rho_b\\&=\rho_b\circ\rho_a\circ\mu_b\circ\mu_a\circ\rho_a\circ\rho_b\\&=\rho_b\circ\rho_a,
	\end{align*}
	and $\mu_a\circ\rho_b=\rho_b\circ\mu_b\circ\mu_a\circ\rho_b=\rho_b\circ\mu_a\circ\mu_b\circ\rho_b=\rho_b\circ\mu_a$.
\end{proof}
For any $P\in\SP(n)$, $A=\{a_1, \ldots, a_m\}\subseteq\MB(P)$ and $B=\{b_1, \ldots, b_s\}\subseteq\Succeq(P)$, we define  $\psi_{A,B}(P)=P'$, where $P'$ is the set partition obtained from $P$ by applying $\mu_a$ for each element $a$ of $A$ and applying $\rho_b$ for each element $b$ of $B$. Thus, $$\psi_{A,B}=\mu_{a_1}\cdots\mu_{a_m}\rho_{b_1}\cdots\rho_{b_s}.$$
By the preceding lemmas there is an equivalence relation in the set $\SP(n)$ defined by two set partitions 
$$P\equiv P' \iff (\exists A\subseteq\MB(P)~ \exists B\subseteq \Succeq(P)\text{ such that } \psi_{A,B}(P)=P').$$
Let $\Gamma[P]$ denote the equivalence class containing $P$. 
\begin{proposition}
	For any $P\in\SP(n)$, we have
	\begin{equation} \label{eqncardGamma}
		|\Gamma[P]|=2^{\mb(P)+\sueq(P)}.
	\end{equation}
	
	Moreover, for any $P'\in\Gamma[P]$, we have \begin{equation}\MB(P')\cup\Succeq(P')=\MB(P)\cup\Succeq(P)\end{equation} and 
	\begin{equation}
		\NMB(P')=\NMB(P).
	\end{equation}
	
\end{proposition}

\begin{proof}
	Equation (\ref{eqncardGamma}) follows directly from the fact that for any $A\subseteq\MB(P)\cup\Succeq(P)$, there exists a unique set partition $P'\in\Gamma[P]$ such that $\MB(P')=A$. The rest follows from Lemma~\ref{lemmaMbtoSuc}.
\end{proof}
\begin{corollary}
	The number of set partitions over $[n]$ having exactly one non-merging block is $2^{n-1}, n\geq1$.
\end{corollary}
\begin{proof}
	The set partition over $[n]$ having exactly one non-merging block and no merging block is the trivial set partition, $12\cdots n$, with $[2,n]$ as set of successions. Thus, $\{P\in\SP(n): \nmb(P)=1\}=\Gamma[12\cdots n]$, and it has indeed size $2^{n-1}$ by~(\ref{eqncardGamma}).
\end{proof}
\subsection{Enumeration results}
In this subsection, we employ the bijections we have defined to give some results on the distribution of $\mb(P)$ and $\sueq(P)$, where $P$ is any set partition over $[n]$.

\bigskip
\begin{lemma}\label{lemmaTRcard}
	For any $A, A', B, B'\subseteq[2,n], n\geq2$ such that $A$ and $B$ are disjoint, and $A'$ and $B'$ are disjoint, and $A\cup B=A'\cup B'$, the cardinalities of the sets $\mathcal{T}_n^A\cap\mathcal{R}_n^B$ and $\mathcal{T}_n^{A'}\cap\mathcal{R}_n^{B'}$ are equal.
\end{lemma}
\begin{proof}
	The map $\psi_{A\backslash A',B\backslash B'}$ yields a bijection between these sets.
\end{proof}
We note that for any disjoint subsets $A$ and $B$ of $[2,n]$, the restriction of $\psi_{A,B}$ to the set $\mathcal{T}_n^A\cap\mathcal{R}_n^B$ provides a bijection between this set and $\mathcal{T}_n^B\cap\mathcal{R}_n^A$. 
(Note that $\psi_{\emptyset,\emptyset}$ restricts to the identity on $\mathcal{T}_n^\emptyset\cap\mathcal{R}_n^\emptyset$). Since the collection of such $\mathcal{T}_n^A\cap\mathcal{R}_n^B$ forms a partition of $\SP(n)$, we can put these restrictions together to obtain an involution $\psi$. In other words, for any set partition $P$, we let $\psi(P)=\psi_{\MB(P),\Succeq(P)}(P)$.
\begin{example}
	Let $P=1~4~5|2~6~7~9|3|8~10$. We have $\MB(P)=\{8\}$, $\Succeq(P)=\{5,7\}$, and $\psi_{\{8\},\{5,7\}}(P)=\rho_5\rho_7\mu_8(P)=P'$. Then $\mu_8(P)=1~4~5|2~6~7~8~10|3~9, ~\rho_7(\mu_8(P))=1~4~5|2~6~9|3|7~8~10, ~ \rho_5(\rho_7(\mu_8(P)))=1~4|2~6~9|3|5|7~8~10{=}P'$ and $\MB(P')=\{5, 7\}, \Succeq(P')=\{8\}, \nmb(P')=3=\nmb(P)$. 
\end{example}
In particular (or by Lemma~\ref{lemmaTRcard}), we have
\begin{theorem}\label{thmmbsucceq}
	Let $n\geq1$ and
	$$F_n(q,t,r)=\sum_{P\in\SP(n)}q^{\mb(P)}t^{\sueq(P)}r^{\nmb(P)}.$$
	Then
	$$F_n(q,t,r)=F_n(t,q,r).$$ \qed
\end{theorem}
\begin{proposition}
	For any $A\subseteq[2,n], n\geq2$, the cardinality of the set $\mathcal{T}_n^A$ is given by $$|\mathcal{T}_n^A|=B(n{-}1{-}|A|).$$
\end{proposition}
\begin{proof}
	Let $P\in\mathcal{T}_n^A$, where $A=\{a_1, \ldots, a_m\}$. If $P'=\mu_{a_1}\cdots\mu_{a_m}(P)$, then $P'\in\mathcal{T}_n^\emptyset$. We then delete each $a\in A$ obtaining a set partition $P''$ on $n{-}|A|$ letters having no merging blocks, i.\ e., $P''\in\mathcal{T}_{n-|A|}^\emptyset$. So the map $P\mapsto P''$ is a bijection, whence, indeed, by Remark \ref{remmerfree} $|\mathcal{T}_n^A|=|\mathcal{T}_{n-|A|}^\emptyset|=B(n{-}1{-}|A|)$.
\end{proof}
Let $\SP^*(n)$ denote the set of all set partitions $P\in\SP(n)$ such that the removal of $n$ creates a new merging block.
\begin{proposition}
	We have
	$$\sum_{\substack{P\in\SP^*(n+2)} }q^{\mb(P)}t^{\sueq(P)}r^{\nmb(P)}=n\sum_{Q\in\SP(n)}q^{\mb(Q)}t^{\sueq(Q)}r^{\nmb(Q)+1}.$$
\end{proposition}
\begin{proof}
	We prove the assertion by providing a bijection between the sets  $[2,n+1]\times\SP(n)$ and $\SP^*(n{+}2)$. Let $\theta: [2,n{+}1]\times\SP(n)\mapsto\SP^*(n{+}2)$ be the map associating $(a,P)$ with the set partition $P'$, where $P'$ is obtained from $(a,P)$ as follows. Increase by $1$ every integer greater than or equal to $a$ in $P$ and insert $a$ into the block containing $a{-}1$. Now apply $\rho_a$ to the resulting set partition, and insert $n{+}2$ in the block preceding the merging block newly created. Note that $P'\in\SP^*(n{+}2)$ and $\theta$ is a bijection such that $\mb(P')=\mb(P), \sueq(P')=\sueq(P)$, and $\nmb(P')=\nmb(P)+1$. Therefore, we have the assertion.
\end{proof}
Let $h_k(n,m,s):=|\{P{\in}\SP(n): \blocks(P){=}k, \mb(P)=m, \sueq(P)=s\}|$. Then $h_1(n,0,n{-}1)=h_n(n,n{-}1,0)=1, n\ge1$, and $h_k(n,m,s)=0$, where $k>n, m\ge k, s\ge n$, or $n<0$. 
\begin{proposition}\label{propmergtonomerg}
	For $n\geq1$, we have
	\begin{equation}\label{mersuctoonlysuc}
		h_k(n,m,s)={{m+s}\choose m}h_{k-m}(n,0,s+m).
	\end{equation}
\end{proposition}
\begin{proof}
	We start with any set partition over $[n]$ having $k{-}m$ blocks and no merging blocks. If the set partition has $m{+}s$ successions, then we can create $m$ merging blocks, by applying the maps $\rho$, in ${{m+s}\choose{m}}$ ways. Thus, by the product rule, we have the result.
\end{proof}
We now give some consequences of the above proposition. 
\begin{proposition}
	Given $n> s\ge1$, we have
	\begin{equation}\label{onesuccfromzero}
		h_k(n,0,s)={n-1\choose s}h_k(n{-s},0,0).
	\end{equation}
\end{proposition} 
\begin{proof}
	Let $P^{(0)}=P\in\mathcal{T}_{n-s}^\emptyset\cap\mathcal{R}_{n-s}^\emptyset$. There are ${n-1\choose s}$ possible ways to choose a subset of $[2,n]$ having size $s$. For any such set $A=\{a_1, \ldots, a_s\}$ with $a_1<\cdots <a_s$ and for $i=1, \ldots, s$, let $P^{(i)}$ be the set partition obtained from $P^{(i-1)}$ by increasing  by $1$ each integer greater than or equal to $a_i$ and  inserting $a_i$ in the  block containing $a_i{-}1$. So $P^{(s)}$ is a set partition over $[n]$ with $\Succeq(P^{(s)})=A$. Hence, by the product rule, we obtain the result.
\end{proof}
By combining (\ref{mersuctoonlysuc}) and (\ref{onesuccfromzero}) we have the following corollary.
\begin{corollary}
	For $n\geq1$ we have
	\begin{equation}\label{eqnmersuctozero}
		h_k(n,m,s)={{n-1}\choose m,s,n{-}m{-}s{-}1}h_{k-m}(n{-}m{-}s,0,0).
	\end{equation}
\end{corollary}
We let $G(x,y,z,w):=\sum\limits_{n,k,m,s\ge0}h_k(n,m,s)x^ny^kz^mw^s$, and $J(x,y):=\sum\limits_{n,k\ge0}h_k(n,0,0)x^ny^k$. Then we have 
\begin{proposition} $G(x,y,z,w)=J(x(1-xyz-xw)^{-1},y)$.
\end{proposition}
\begin{proof} 
	By (\ref{eqnmersuctozero}) we have indeed
	\begin{align*}
		G(x,y,z,w)&=\sum\limits_{n,k,m,s\ge0}{n-1+m+s\choose m,s,n-1}x^{m+s}y^mz^mw^sh_k(n,0,0)x^ny^k\\
		&=\sum\limits_{n,k\ge0}\sum\limits_{m\ge0}\sum\limits_{s\ge0}{n-1+m+s\choose n-1}{m+s\choose m}(xyz)^{m+s}y^mz^mw^sh_k(n,0,0)x^ny^k\\
		&=\sum\limits_{n,k\ge0}\sum\limits_{m+s\ge0}{n-1+m+s\choose n-1}\sum\limits_{m\ge0}{m+s\choose m}(xyz)^{m+s}y^mz^mw^sh_k(n,0,0)x^ny^k\\
		&=\sum\limits_{n,k\ge0}\sum\limits_{m+s\ge0}{n-1+m+s\choose n-1}(xyz+xw)^{m+s}h_k(n,0,0)x^ny^k\\
		&=\sum\limits_{n,k\ge0}\frac{1}{(1-xyz-xw)^n}h_k(n,0,0)x^ny^k\\
		&=J\left(\frac{x}{1-xyz-xw},y\right).
	\end{align*}
\end{proof}
We let $\SP^0(n):=\mathcal{T}_n^\emptyset\cap\mathcal{R}_n^\emptyset$, the set of set partitions having no merging blocks and no successions. So we have $|\SP^0(n,k)|=h_k(n,0,0)$. 
\begin{theorem}
	The number $h_k(n,0,0)$ satisfy the following recurrence relation for all positive integers $n, k$, $n\geq2, 1\leq k\leq \lfloor \frac{n-1}{2}\rfloor$:
	\begin{equation}\label{eqnnomergnosucc}
		h_k(n,0,0)=(k-1)h_k(n-1,0,0)+(n-2)h_{k-1}(n-2,0,0),
	\end{equation}
	where $h_0(n,0,0)=\delta_{n,0}, ~h_1(1,0,0)=1$.
\end{theorem}
\begin{proof}
	Let $n$ and $k$ be fixed positive integers. Let $\SP^0(n,k)=\mathcal{M}\cup\mathcal{N}$, where $\mathcal{M}$ is the subset of $\SP^0(n,k)$ consisting of those set partitions whose removal of $n$ does not create a merging block, and $\mathcal{N}=\SP^0(n,k)\backslash\mathcal{M}$.
	
	Let $P\in\SP^0(n{-}1,k)$ and $P'$ be the set partition obtained from $P$ by inserting $n$ to any of its blocks except the block containing $n{-}1$. Then $P'\in\mathcal{M}$. Since there are $k{-}1$ possibilities where to insert $n$, we have the first term of the right-hand side of (\ref{eqnnomergnosucc}). 
	
	On the other hand, consider $a\in[2,n{-}1]$ and $P\in\SP^0(n{-}2,k{-}1)$. Let $\kappa$ be the map that associates $(a,P)$ with the set partition $P'$ obtained as follows. Increase all integers greater than or equal to $a$ in $P$ by $1$, split the rightmost block containing element(s) of the set $[a{-}1]$ after the rightmost element of $[a{-}1]$, insert $n$ and $a$ to the left and the right blocks of the splitted block, respectively. Then let the resulting partition be $P^*$. If $a+1$ is a succession in $P^*$, then let $P'= \Swap_{a+1}^{(i,j)}(P^*)$, where $(i,j)$ is the pair of indices of the blocks containing $a$ and $a{-}1$ in $P^*$; Otherwise, let $P'=P^*$. It can then be seen that $\kappa: [2,n{-}1]\times\SP^0(n{-}2,k{-}1)\mapsto\mathcal{N}$ is a bijection. Therefore, we have $|\mathcal{N}|=(n{-}2)h_{k-1}(n{-}2,0,0)$, the second term of the right-hand side of (\ref{eqnnomergnosucc}).
\end{proof}
Up to a shift on both $n$ and $k$, this is the same sequence as OEIS entry number \seqnum{A008299}, counting set partitions without singletons. Therefore, there should be a natural bijection between these sets though so far we couldn't find one.

We now consider the distribution of the number of successions in a set of merging-free partitions $\mathcal{T}_n^\emptyset$ having a fixed number of blocks.
\begin{theorem}
	The numbers $h_k(n,0,s)$ satisfy the following recurrence relation for all positive integers $n, k, s$, $1\leq s\leq n-2k+1, 1\leq 2k-1\leq n$:
	\begin{equation}\label{rronlysucc}
		h_k(n,0,s)=h_k(n{-}1,0,{s{-}1})+(k{-}1)h_k(n{-}1,0,s)+(s{+}1)h_{k-1}(n-1,0,{s+1});
	\end{equation} 
	and $h_k(n,0,0)$ satisfies (\ref{eqnnomergnosucc}). 
\end{theorem}
\begin{proof}
	It is possible to obtain any set partition $P'\in\mathcal{T}_n^\emptyset$ recursively either from $P\in\mathcal{T}_{n-1}^\emptyset$ by inserting $n$ in any of the existing blocks of $P$, or from any $P^*\in\mathcal{T}_{n-1}^{\{a\}}$, where $a\in[2,n]$, by inserting $n$ in the block preceding a merging block of $P^*$. In the first case, if $n$ is inserted into the block containing $n{-}1$, then the number of successions increases by $1$, but otherwise it remains the same; anyhow the number of blocks remains the same. This explains the first two terms of the right-hand side of (\ref{rronlysucc}). In the second case, $P:=\mu_a(P^*)$ has $\Succeq(P)=\Succeq(P')\cup\{a\}$ and the number of blocks one less than that of $P'$. Since $P^*=\rho_a(P)$ and $a$ has $\sueq(P)$ possibilities, this yields the third term. 
\end{proof}
\begin{proposition} 
	Let $H_k(x,z)=\sum\limits_{n\ge 2k-1}\sum\limits_{s\geq0}h_k(n,0,s)z^sx^n$. Then we have
	\begin{equation}\label{eqnHk}
		H_k(x,z)=\frac{x}{1-x(k{-}1+z)}\frac{\partial}{\partial z}\left(H_{k-1}(x,z)\right), k\ge2.
	\end{equation}
\end{proposition}
\begin{proof} 
	We define the polynomial $H(n,k;z)=\sum\limits_{s=0}^{n-1}h_k(n,0,s)z^s$. Then by (\ref{rronlysucc}) we have
	\begin{align*}
		\sum_{s\geq1}h_k(n,0,s)z^s&=\sum_{s\geq1}h_k(n{-}1,0,s{-}1)z^s+\sum_{s\geq1}(s{+}1)h_{k-1}(n{-}1,0,s{+}1)z^s\\&~~~+\sum_{s\geq1}(k{-}1)h_k(n{-}1,0,s)z^s\\
		H(n,k;z)-h_k(n,0,0)&=zH(n{-}1,k;z)+H_z(n{-}1,k{-}1;z)-h_{k-1}(n{-}1,0,1)\\&~~~+(k{-}1)(H(n{-}1,k;z)-h_k(n{-}1,0,0)).
	\end{align*}
	By applying (\ref{eqnnomergnosucc}) and (\ref{onesuccfromzero}) we obtain
	\begin{equation}\label{eqnwithoutS}
		H(n,k;z)=(z{+}k{-}1)H(n{-}1,k;z)+\frac{\partial}{\partial z}H(n{-}1,k{-}1;z), n\ge 2k-1.
	\end{equation}
	We now let $H_k(x,y)=\sum\limits_{n\ge 2k-1}H(n,k;z)x^n$. Then multiplying (\ref{eqnwithoutS}) by $x^n$ and taking the sum over all $n\ge 2k{-}1$, we obtain  
	\begin{equation*}
		H_k(x,z)=\frac{x}{1-x(k{-}1+z)}\frac{\partial}{\partial z}(H_{k-1}(x,z)).
	\end{equation*}
\end{proof} 
We introduce the following definition.
\begin{definition} \label{vectordefn}
	Let $r\ge0$ and $v=(v_0,v_1,\ldots, v_r)$ be a vector of non-negative integers such that $\sum\limits_{j=0}^rv_j=r$ and for $1\le i\le r$, $s_i(v)>i-2$, where $s_i(v):=\sum\limits_{j=0}^{i-1}v_j$.
\end{definition}
For any such vector $v=(v_0,v_1,\ldots, v_r)$, $v_r\le1$. If $v_r=1$, then we let $v^{(r)}=(v_0,v_1,\ldots, v_{r-1})$. If $v_r=0$, then for $0\le t\le r{-}1$ and $v_t>0$, let $v^{(t)}=(v_0',v_1',\ldots, v_{r-1}')$ be the vector obtained from $v$ by setting $v_t'=v_t{-}\delta_{t,q}$, and deleting $v_r$. We also let $P_v:=\prod_{i=1}^r(s_i(v){-}i+2)$, correspondingly for $P_{v^{(t)}}$. 

We give the following lemma that will be used to prove Theorem \ref{thmSolutionH}.
\begin{lemma}\label{lemmaPv}
	For $r\ge0$ we have 
	\begin{equation*} P_v=\begin{cases} P_{v^{(r)}}, \text{ if } v_{r}=1;\\ \sum\limits_{j=0}^{r-1}(v_j+\delta_{j,0})P_{v^{(j)}}, \text{ if } v_{r}=0.
		\end{cases}
	\end{equation*}
\end{lemma}
\begin{proof}
	If $v_{r}=1$, then the last factor of $P_v=\prod_{i=1}^{r}(s_i(v){-}i+2)$ is $(v_0+\cdots+ v_{r-1}-r+2)=r-v_{r}-r+2=1$. Therefore, $P_v=\prod_{i=1}^{r-1}(s_i(v){-}i+2)=P_{v^{(r)}}$. 
	
	We now assume that $v_{r}=0$.	By definition $P_{v^{(t)}}=\prod_{i=1}^{r-1}(s_i(v^{(t)}){-}i+2)$ and $$v^{(t)}=(v_0,v_1,\ldots, v_t{-}1,v_{t+1},\ldots, v_{r-1}).$$ Then 
	\begin{align}\label{eqnPvt}
		P_{v^{(t)}}&=\prod_{i=1}^{r-1}(s_i(v^{(t)}){-}i+2)\notag \\
		&=\prod_{i=1}^{t}(s_i(v){-}i+2)\cdot\prod_{i=t+1}^{r-1}(s_i(v){-}i+1) 
	\end{align}
	We first use induction on $t$ to prove that \begin{align}\label{eqnsumtoT}
		\sum\limits_{j=0}^{t}(v_j+\delta_{j,0})P_{v^{(j)}}&=\prod_{i=1}^{t+1}(s_i(v){-}i+2)\cdot\prod_{i=t+1}^{r-1}(s_i(v){-}i+1).
	\end{align}
	Observe that 
	$$\sum\limits_{j=0}^{0}(v_j+\delta_{j,0})P_{v^{(j)}}=\prod_{i=1}^{1}(s_i(v){-}i+2)\cdot\prod_{i=1}^{r-1}(s_i(v){-}i+1),$$ 
	and the assertion is true for $t=0$. Suppose that $t>0$, and 
	\begin{align*}\sum\limits_{j=0}^{t-1}(v_j+\delta_{j,0})P_{v^{(j)}}&=\prod_{i=1}^{t}(s_i(v){-}i+2)\cdot\prod_{i=t}^{r-1}(s_i(v){-}i+1).
	\end{align*}
	Now by the induction assumption and (\ref{eqnPvt}) we have
	\begin{align*} \sum\limits_{j=0}^{t}(v_j+\delta_{j,0})P_{v^{(j)}}&=\sum\limits_{j=0}^{t-1}(v_j+\delta_{j,0})P_{v^{(j)}}+v_tP_{v^{(t)}}\\&=\prod_{i=1}^{t}(s_i(v){-}i+2)\cdot\prod_{i=t}^{r-1}(s_i(v){-}i+1)+v_t\left(\prod_{i=1}^{t}(s_i(v){-}i+2)\cdot\prod_{i=t+1}^{r-1}(s_i(v){-}i+1)\right)\\
		&=\prod_{i=1}^{t}(s_i(v){-}i+2)(s_t(v){-}t+1+v_t)\prod_{i=t+1}^{r-1}(s_i(v){-}i+1)\\
		&=\prod_{i=1}^{t+1}(s_i(v){-}i+2)\prod_{i=t+1}^{r-1}(s_i(v){-}i+1).
	\end{align*}
	and thus (\ref{eqnsumtoT}) is proved. Then (\ref{eqnsumtoT}) for $t=r{-}1$ and the definition of $P_v$ yields the result of the lemma for which $v_{r}=0$.
\end{proof}
\begin{theorem}\label{thmSolutionH}
	The generating function for $H_k(x,z)$ is given by
	\begin{equation}\label{eqnHkSolution}
		H_k(x,z)=\frac{x^{2k-1}}{(1-xz)\Pi_{j=0}^{k-1}(1-x(j+z))}\sum_v\frac{\Pi_{i=1}^{k-2}(s_i(v){-}i+2)}{\Pi_{j=0}^{k-2}(1-x(j+z))^{v_j}}, k\ge2,
	\end{equation}
	where $v=(v_0,v_1,\ldots, v_{k-2})$ as in Definition \ref{vectordefn} (with $r=k{-}2$) with $H_0(x,z)=1$ and $H_1(x,z)=\frac{x}{1-xz}$.
\end{theorem}
\begin{proof}
	For $0\le j\le k-2$, let $a_j:=1-x(j+z)$ and $a^v:=a_0^{v_0}\cdots a_{k-2}^{v_{k-2}}$. Then from (\ref{eqnHk}) we have $H_k(x,z)=\frac{x}{a_{k-1}}\frac{\partial}{\partial z}(H_{k-1}(x,z))$, and the right-hand side of (\ref{eqnHkSolution}) is 
	\begin{equation*} \frac{x^{2k-1}}{a_0^2a_1\cdots a_{k-1}}\sum\limits_{v}\frac{P_v}{a^v}=\frac{x^{2k-1}}{a_0^2a_1\cdots a_{k-1}}\left(\sum\limits_{\substack{v,\\ v_{k-2}=1}}\frac{P_v}{a^v}+\sum\limits_{\substack{v,\\ v_{k-2}=0}}\frac{P_v}{a^v}\right).
	\end{equation*}
	By applying Lemma \ref{lemmaPv}, we have
	$\sum\limits_{\substack{v,\\ v_{k-2}=1}}\frac{P_v}{a^v}=\sum\limits_{\substack{v,\\ v_{k-2}=1}}\frac{P_{v^{(k-2)}}}{a^v}$, and 
	\begin{align*}\sum\limits_{\substack{v,\\ v_{k-2}=0}}\frac{P_v}{a^v}&=\sum\limits_{\substack{v,\\ v_{k-2}=0}}\left(\frac{(v_0{+}1)P_{v^{(0)}}{+}v_1P_{v^{(1)}}{+}\cdots{+}v_{k-2}P_{v^{(k-2)}}}{a^v}\right)\\
		&=\sum\limits_{\substack{v,\\ v_{k-2}=0}}\left(\frac{2P_{v^{(0)}}{+}P_{v^{(1)}}{+}\cdots{+}P_{v^{(k-2)}}}{a^v}+\frac{(v_0{-}1)P_{v^{(0)}}{+}(v_1{-}1)P_{v^{(1)}}{+}\cdots{+}(v_{k-2}{-}1)P_{v^{(k-2)}}}{a^v}\right)\\
		&=\sum\limits_{\substack{v,\\ v_{k-2}=0}}\left(\frac{2}{a_0}\frac{P_{v^{(0)}}}{a^{v^{(0)}}}+\frac{1}{a_1}\frac{P_{v^{(1)}}}{a^{v^{(1)}}}+\cdots+\frac{1}{a_{k-2}}\frac{P_{v^{(k-2)}}}{a^{v^{(k-2)}}}+\left(\frac{v_0{-}1}{a_0}\frac{P_{v^{(0)}}}{a^{v^{(0)}}}+\cdots+\frac{v_{k-3}{-}1}{a_{k-3}}\frac{P_{v^{(k-3)}}}{a^{v^{(k-3)}}}\right)\right)
	\end{align*}
	Therefore, 
	\begin{align*}\frac{x^{2k-1}}{a_0^2\Pi_{j=1}^{k-1}a_j}\sum\limits_{v}\frac{P_v}{a^v}&=\frac{x^{2k-1}}{a_0^2\Pi_{j=1}^{k-1}a_j}\left(\left(\frac{2}{a_0}{+}\frac{1}{a_1}{+}\cdots{+}\frac{1}{a_{k-2}}\right)\sum\limits_{v^{(t)}}\frac{P_{v^{(t)}}}{a^{v^{(t)}}}{+}\sum\limits_{v^{(t)}}\frac{P_{v^{(t)}}}{a^{v^{(t)}}}\left(\frac{v_0{-}1}{a_0}{+}\cdots{+}\frac{v_{k-3}{-}1}{a_{k-3}}\right)\right) \\
		&=\frac{x}{a_{k-1}}\frac{\partial}{\partial z}\left(\frac{x^{2k-3}}{a_0^2a_1\cdots a_{k-2}}\sum\limits_{v^{(t)}}\frac{P_{v^{(t)}}}{a^{v^{(t)}}}\right)\\
		&=\frac{x}{a_{k-1}}\frac{\partial}{\partial z}(H_{k-1}(x,z))\\
		&=H_k(x,z),
	\end{align*}
	indeed.
\end{proof}
By the fact that $\sum\limits_{n\ge k}S(n,k)x^n=\frac{x^k}{\Pi_{j=0}^{k}(1-jx)}, k\ge0$, we have
\begin{corollary}
	\begin{equation*}
		H_k(x,0)=x^k\sum\limits_{n\ge k-1}S(n,k-1)x^n\sum_v\frac{\Pi_{i=1}^{k-2}(s_i(v){-}i+2)}{\Pi_{j=0}^{k-2}(1-jx)^{v_j}}, k\ge2,
	\end{equation*}
	with $H_0(x,0)=1, H_1(x,0)=x$.
\end{corollary}
Let us use the notation $h_{n,m,s}:=\sum_{k=1}^nh_k(n,m,s)$.
\begin{proposition}
	For $n\geq1$ we have
	\begin{equation*}
		\sum_{s=0}^{n-1}2^sh_{n,0,s}=B(n),
	\end{equation*}
	where $B(n)$ is the $n$'th Bell number.
\end{proposition}
\begin{proof}
	By Proposition \ref{propmergtonomerg} we have $h_{n,m,s-m}={s\choose m}h_{n,0,s}$. Thus, $$\sum_{m=0}^sh_{n,m,s-m}=\sum_{m=0}^s{s\choose m}h_{n,0,s}=2^sh_{n,0,s}.$$ 
	Hence taking the sum over all possible $s$ we have the result.
\end{proof}		
\section{Acknowledgements}
The first author is grateful for the financial support extended by the cooperation agreement between the International Science Program (ISP) at Uppsala University and Addis Ababa University, the support by the CDC-Simons for Africa, and IRIF. We appreciate the hospitality we got from Stockholm University during the research visit of the first author. We also thank our colleagues from CoRS (Combinatorial Research Studio) for valuable discussions and comments, in particular, we thank Dr.\ Per Alexandersson of Stockholm University for his crucial discussions and suggestions.
%%%%%%%%%%%%%%%%%%%%%%%%%%%%%%%%%%%%%%%	


\begin{thebibliography}{99}
	\addcontentsline{toc}{chapter}{References}
	\bibitem{Ba1} J. L. Baril, \emph{Gray code for permutations with fixed number of cycles}, Discrete Mathematics, III (2006), III-III.
	\bibitem{Ba} J. L. Baril, \emph{Statistics-preserving bijections between classical and cyclic permutations}, Information Processing Letters, 113 (2013), 17-22.
	\bibitem{Ba-Va}J. Baril and V. Vajnovszki, \textit{A permutation code preserving a double Eulerian bistatistic}, Discrete Applied Mathematics, 224 (2017) 9-15, (2018).
	\bibitem{Be-Ma} F. Beyene and R. Mantaci, \textit{Merging-Free Partitions and Run-Sorted Permutations}, \url{http://arxiv.org/abs/2101.07081}.
	\bibitem{Be-Ma2} F. Beyene and R. Mantaci, \textit{Permutations with non-decreasing transposition array and pattern avoidance}, \url{http://arxiv.org/abs/2111.11527}.
	\bibitem{Bo} M. Bona, \textit{Introduction to Enumerative Combinatorics}, The McGraw Hill Companies, (2007).
	\bibitem{Ca} D. Callan, On conjugates for set partitions and integer compositions, \url{arxiv:math/0508052v3 [math.CO] 11 Oct 2005}.
	\bibitem{Du-Vi} D. Dumont and G. Viennot, \textit{A combinatorial interpretation of the Seidel generation of Genocchi numbers}, Ann.
	Discrete Math. 6, 77-87, (1980).
	\bibitem{Fo-Ze} D. Foata and D. Zeilberger, \textit{Denert's Permutation Statistic is indeed Euler-Mahonian}, Studies in Applied Mathematics, 31-59, (1990).
	\bibitem{Le}D. H. Lehmer, \textit{Teaching combinatorial tricks to a computer}, in Proc. Sympos. Appl. Math., 10 (1960), Amer. Math. Soc., 179-193.
	\bibitem{Ma-Rak} R. Mantaci and F. Rakotondrajao. A permutation representation that knows what ``Eulerian'' means. Discrete Mathematics and Theoretical Computer Science 4, 2001, 101-108.
	\bibitem{Ma} T. Mansour. Combinatorics of set partitions. Taylor \& Francis Group, LLC, 2013.
	\bibitem{Ma-Mu} T. Mansour and A. O. Munagi, Set partitions with circular successions, European Journal of Combinatorics, 42 (2014), 207-216.
	\bibitem{Ma-Ra} T. Mansour and R. Rastegar, Fixed points of a random restricted growth sequence, \url{arxiv:2012.06891v2 [math.CO] 24 Jun 2021}.
	\bibitem{Mu} A. O. Munagi, Set Partitions with successions and separations, Intl. J. Math. Math. Sci. (2005), 451-463.
	
	\bibitem{Po-Va} M. Poneti and V. Vajnovszki, \textit{Generating restricted classes of involutions, Bell and Stirling permutations}, European Journal of Combinatorics 31, 553-564, (2010).
	\bibitem{Ro} G. Rota, \textit{The number of partitions of a set}, Amer. Math. Monthly 71, 498–504, (1964).
	\bibitem{Sl} N. J. A. Sloane et al., The on-line encyclopedia of integer sequences, Available at https://oeis.org, 2020
	\bibitem{St1} R. P. Stanley, \textit{Enumerative Combinatorics}, Vol. 1, $2^{ed}$ ed, Cambridge Studies of Advanced Mathematics, Cambridge University Press, (2011).
\end{thebibliography}
\end{document}